\documentclass[12pt]{article} 
\usepackage{amsmath}
\usepackage{amssymb} 
\usepackage{amscd}
\usepackage{amsthm}
\usepackage{mathabx}
\usepackage[title,titletoc,toc]{appendix}
\usepackage[utf8]{inputenc}
\usepackage{stmaryrd}
\usepackage{hyperref}
\usepackage{fullpage}

\newtheorem{theorem}{Theorem}

\newtheorem{proposition}[theorem]{Proposition}
\newtheorem{lemma}[theorem]{Lemma}
\newtheorem{corollary}[theorem]{Corollary}

\def\ad{\mathrm{ad}}
\def\gr{\mathrm{gr}}

\def\dim{{\mbox{dim}}}

\def\ker{{\mbox{Ker}}}     
\def\Der{{\mbox{Der}}}

\def\Hom {{\mbox{Hom}}}

\def\End{{\mbox{End}}}

\def\Lie{{\mbox{Lie}}}

\def\cali{{\cal I}}
 
\def\cala{{\cal A}}

\def\calm{{\cal M}}

\def\calh{{\cal H}}

 \def\fracg{{\mathfrak g}}

\def\fracT{{\mathfrak {T}}}

\def\bbbone{\mbox{\rm 1\hspace {-.6em} l}}

\numberwithin{equation}{section}

\begin{document}
\enlargethispage{3cm}

\thispagestyle{empty}
\begin{center}
{\bf THE WEIL ALGEBRA OF A HOPF ALGEBRA}\\
\vspace{0,5cm}
{\bf \large I - A noncommutative framework}
\end{center}

\vspace{0.3cm}

\begin{center}
Michel DUBOIS-VIOLETTE
\footnote{Laboratoire de Physique Th\'eorique, UMR 8627\\
Universit\'e Paris XI,
B\^atiment 210\\ F-91 405 Orsay Cedex\\
Michel.Dubois-Violette$@$u-psud.fr} and 
Giovanni LANDI \footnote{Dipartimento di Matematica, Universit\`a di Trieste\\
 Via A. Valerio 12/1, I-34127 Trieste, Italy\\
  and INFN, Sezione di Trieste, Trieste, Italy\\
landi$@$univ.trieste.it \\
Partially supported by the project {\sl PRIN08 - Noncommutative Geometry, Quantum Groups and Applications}.

 \vspace{1,5cm}} 
\end{center}
\vspace{0,5cm}
 \vspace{0,5cm}

\begin{abstract}
We generalize the notion, introduced by Henri Cartan, of an operation of a Lie algebra $\fracg$ in a graded differential algebra $\Omega$. We define the notion of an operation of a Hopf algebra $\calh$ in a graded differential algebra $\Omega$ which is refered to as a $\calh$-operation. We then generalize for such an operation the notion of algebraic connection.
Finally we discuss the corresponding noncommutative version of the Weil algebra: The Weil algebra $W(\calh)$ of the Hopf algebra $\calh$ is the universal initial object of the category of $\calh$-operations with connections.
\end{abstract}

\vfill
\today \\
\\
\newpage
\enlargethispage{\baselineskip}
\tableofcontents
\parskip 1ex

\newpage
\section{Introduction}
 This paper is devoted to the noncommutative version of the notion of Cartan operation of a Lie algebra $\fracg$ in a graded differential algebra $\Omega$, 
 \cite{car:1951a}, \cite{car:1951b}. Our aim is to generalize the classical theory to the noncommutative setting in two respects: firstly to replace the Lie algebra $\fracg$, or more precisely its universal enveloping algebra $U(\fracg)$, by an arbitrary Hopf algebra $\calh$, and secondly to allow arbitrary (noncommutative) graded differential algebras.

We define the notion of an operation of a Hopf algebra $\calh$ in a graded differential algebra $\Omega$. By letting $\calh$ be the universal enveloping algebra $U(\fracg)$ of a Lie algebra $\fracg$ we obtain by restriction to $\fracg$ an operation of $\fracg$ in $\Omega$ as defined by Henri Cartan (except that $\Omega$ is not assumed graded commutative). We discuss the converse direction, namely from an operation of $\fracg$ to an operation of $U(\fracg)$.
We observe that various extensions of the Cartan calculus for the differential calculi on quantum groups (see e.g. \cite{sch-wat:1994}, 
\cite{kli-sch:1997} and \cite{asc-cas:1993}) fall in this framework. We then introduce the notion of an algebraic connection for an operation of a Hopf algebra $\calh$ in a graded differential algebra $\Omega$ generalizing thereby the corresponding notion induced by Henri Cartan \cite{car:1951a} (see also in \cite{kam-ton:1975}, \cite{gre-hal-van:1976}) for an operation of a Lie algebra in a graded \underbar{commutative} differential algebra.
The category of operations with connections of a given Hopf algebra $\calh$ in graded differential algebras has a universal initial object $W(\calh)$ which we describe and which is the appropriate generalization of the Weil algebra. 

We stress that our $W(\calh)$ is not connected with the noncommutative Weil algebra of \cite{ale-mei:1999}, even when $\calh=U(\fracg)$. In fact, the purpose of the nice construction of \cite{ale-mei:1999} is different from that of this paper.
Our construction also differ from the one in the interesting paper \cite{cra:2002}.
The relation of our construction with the one of \cite{cra:2002} is discussed at the end of this paper.

We shall start with a summary of the classical theory. In this summary we deliberately drop two assumptions in the definition of an operation of a Lie algebra $\fracg$ in a graded differential algebra $\Omega$. The first one is the axiom $(i_X)^2=0$, for any $X\in \fracg$, and the second one is the graded commutativity of $\Omega$. Indeed, as explained in the conclusion $(i_X)^2=0$, 
for any $X\in \fracg$, follows from the other axioms in all cases of interest and plays no role otherwise while the operation of $\fracg$ in a noncommutative $\Omega$ makes sense and is useful as pointed out for instance in \cite{mdv:1988} for the case where $\fracg$ is the Lie algebra $\Der(\cala)$ of all derivations of an algebra $\cala$ and where $\Omega$ is the universal differential calculus $\Omega(\cala)$ over $\cala$, \cite{ac:1982}, \cite{ac:1994},\cite{kar:1983}, \cite{kar:1987}, \cite{mdv:2001}.

As mentioned in \S\ref{defope}, the standard example comes from the theory of differential forms on principal bundles. 
Then, the analogues of this classical example are the appropriate differential calculi over noncommutative principal bundles. Among these noncommutative principal bundles let us mention the principal bundles in \cite{lan-sui:2008} on the noncommutative manifolds of \cite{ac-lan:2001}, \cite{ac-mdv:2002a}: in particular the $SU(2)$ principal bundle \cite{lan-sui:2005} over a noncommutative 4-sphere. The $SU_q(2)$-principal bundle over a quantum 4-sphere in \cite{bra-lan:2011} uses in a crucial way the covariant calculi of \cite{wor:1987b}, \cite{wor:1989} and generalizes the $U(1)$-fibration \cite{brz-maj:1993}, \cite{brz-maj:1995} over the quantum 2-sphere. It is also worth noticing here that the theory of Hopf-Galois extensions includes a general formulation of noncommutative principal bundles, see e.g. \cite{sch:2004} and references therein.

Throughout this paper $\mathbb K$ denotes a field and all vector spaces and algebras are over $\mathbb K$. By an algebra (resp. a Lie algebra) without other specification we always mean a unital associative algebra (resp. a finite-dimensional Lie algebra); the unit of such an algebra will be denoted by $\bbbone$ whenever no confusion arises. Except in the appendix where $\mathbb Z$-graduations are considered, 
by a graded algebra, we mean a $\mathbb N$-graded algebra $\cala=\oplus_{n\geq 0} \cala_n$. Given a vector space $E$, its dual is denoted by $E^\ast$ and given a linear mapping $\varphi:E\rightarrow F$ we denote by $\varphi^t:F^\ast\rightarrow E^\ast$ the corresponding transposed linear mapping. The tensor algebra of $E$ is denoted by $T(E)=\oplus_n E^{\otimes^n}$, the symmetric algebra of $E$ is denoted by $SE=\oplus_n S^nE$ and the exterior algebra of $E$ is denoted by $\wedge E=\oplus_n{\wedge^nE}$. The Koszul rule is understood for tensor products of linear mappings between graded vector spaces and we use Sweedler notation for coproducts:  $\Delta h=\sum_j h^{(1)}_j \otimes h^{(2)}_j$, with $^{(1)}$ and $^{(2)}$ denoting respectively the first and the second leg components in the tensor product.

 \section{Operations of Lie algebras}\label{cathe}

In this section we review basic definitions and facts on the notion of operation introduced by Henri Cartan in 1950, \cite{car:1951a}, \cite{car:1951b}, (see also \cite{gre-hal-van:1976}, \cite{mdv:2001}) and describe some related developments. 

 \subsection{Definition of $\fracg$-operations}  \label{defope}
 Let $\fracg$ be a Lie algebra and let $\Omega$ be a graded differential algebra with differential denoted $d$.
 
 An {\sl operation of the Lie algebra $\fracg$ in the graded differential algebra} $\Omega$ is a linear mapping
 \begin{equation}
 X\mapsto i_X
 \label{inder}
 \end{equation}
 of $\fracg$ into the space $\Der^{-1}(\Omega)$ of the antiderivations (graded-derivations) of degree -1 of $\Omega$ such that if one defines the derivation $L_X$ of degree 0  by
 \begin{equation}
 L_X=i_Xd+di_X
 \label{Lieder}
 \end{equation}
 for $X\in \fracg$, then one has
 \begin{equation}
 [i_X,L_Y]=i_{[X,Y]}
 \label{axop}
 \end{equation}
 for any $X,Y\in \fracg$. It follows from (\ref{Lieder}) and (\ref{axop}) that one has
 \begin{equation}
 [L_X,d]=0
 \label{dhom}
 \end{equation}
 and
 \begin{equation}
 [L_X,L_Y]=L_{[X,Y]}
 \label{Liehom}
 \end{equation}
 for any $X,Y\in \fracg$. Relation (\ref{Liehom}) means that one has a homomorphism of Lie algebras $L$ of $\fracg$ into the Lie algebra $\Der^0(\Omega)$ of all derivations of degree 0 of $\Omega$.
 
 An element $\alpha$ of $\Omega$ is said to be {\sl invariant} if one has 
 $L_X(\alpha)=0$
 for any $X\in \fracg$ while $\alpha$ is said to be {\sl horizontal} if one has
$i_X(\alpha)=0$
 for any $X\in \fracg$. Finally, $\alpha\in \Omega$ is said to be {\sl basic} if it is both invariant and horizontal.  
  
The set $\Omega_I$ of all invariant elements of $\Omega$ is a graded differential subalgebra of $\Omega$, the set $\Omega_B$ of all basic elements of $\Omega$ is a graded differential subalgebra of $\Omega_I$ (and therefore of $\Omega$) while the set $\Omega_H$ of all horizontal elements of $\Omega$ is only a graded subalgebra of $\Omega$ which is stable by the $L_X$ $(X\in \fracg)$. The cohomology $H_I(\Omega)$ of $\Omega_I$ and the cohomology $H_B(\Omega)$ of $\Omega_B$ are refered to respectively as the {\sl invariant cohomology} and the {\sl basic cohomology} of $\Omega$ (whenever no confusion arises concerning the operation).

Of course, the whole terminology above comes from the theory of differential forms on principal bundles. Let $G$ be a Lie group with Lie algebra $\mathrm{Lie}(G)=\fracg$ and let $P=P(M,G)$ be a principal $G$-bundle over $M$ (the basis) \cite{kob-nom:1963} with projection $\pi:P\rightarrow M$. The projection $\pi$ induces a projection $T(\pi):T(P)\rightarrow T(M)$ of the tangent bundle of $P$ onto the tangent bundle of $M$ and induces by duality an injective homomorphism $\pi^\ast:\Omega(M)\rightarrow \Omega(P)$ of graded differential algebras  of the space $\Omega(M)$ of differential forms on the basis $M$ into the space $\Omega(P)$ of differential forms on $P$. The image $\pi^\ast(\Omega(M))$  of $\pi^\ast$ is denoted by $\Omega_B(P)$ and its elements are called basic differential forms on $P$. A tangent vector to $P$ is said to be vertical whenever its projection on $T(M)$ via $T(\pi)$ vanishes. To each $X\in \fracg$ corresponds a fundamental vector field on $P$ which is vertical and which we also denote by $X$. The horizontal forms on $P$ are the forms on $P$ such that the inner derivations with the vertical vector fields vanish. Finally a form on $P$ invariant by the action of the structure group is said to be invariant. Let $i_X$ denote the inner derivation of $\Omega(P)$ by the fundamental vector field corresponding to $X\in\fracg$. Then one verifies that $X\mapsto i_X$ defines an operation of $\fracg$ in $\Omega(P)$ and that the notions of basicity, horizontality and invariance for the elements of $\Omega(P)$ correspond to the ones associated to this operation.  

It is worth noticing here that the notion of Cartan operation appears in some related example, for instance it plays a fundamental role in the computation of the local BRS cohomology of gauge theory, \cite{mdv-tal-via:1985b}, \cite{mdv:1987b}, \cite{mdv-hen-tal:1992}.

A graded differential algebra $\Omega$ equipped with an operation of the Lie algebra $\fracg$ as above will be refered to as a {\sl $\fracg$-operation}.  There is an obvious notion of morphism for $\fracg$-operations. This defines {\sl the category of $\fracg$-operations}.

\subsection{$\fracg$-operations in the differential envelopes}\label{univop}

Given an algebra $\cala$ we denote by $\Omega(\cala)$ the universal differential calculus over $\cala$. This graded differential algebra $\Omega(\cala)$ is also refered to as {\sl the differential envelope of $\cala$}.

The following lemma is easy to prove by using the universal properties of $\Omega(\cala)$ and of the derivation $d:\cala\rightarrow \Omega^1(\cala)$, \cite{mdv:1988}.

\begin{lemma}\label{UOp}
Let $X\mapsto L^{(0)}_X$ be a homomorphism of Lie algebra of $\fracg$ into the Lie algebra $\Der(\cala)$ of all derivations of an algebra $\cala$ into itself. Then there is a unique operation of $\fracg$ in $\Omega(\cala)$ such that $L_X(a)=L^{(0)}_X(a)$ for any $a\in \cala$.
\end{lemma}

As a corollary of this lemma, one has the following theorem.

\begin{theorem}\label{eOp}
Let $\Omega$ be generated in degree $\mathrm{0}$ by $\cala$ as graded differential algebra. Then any operation of $\fracg$ in $\Omega$ is the quotient of a unique operation of $\fracg$ in $\Omega(\cala)$.
\end{theorem}

In other words any operation of $\fracg$ in a graded differential algebra generated in degree 0 can be identified with an operation of $\fracg$ in a differential envelope (through the corresponding canonical surjective homomorphism of graded differential algebras).

\subsection{From $L_X$ for $X\in \fracg$ to $L_h$ for $h\in U(\fracg)$}\label{extL}

Let $\Omega$ be a $\fracg$-operation. The linear mapping $X\mapsto L_X$ is, in view of (\ref{Liehom}), a representation of the Lie algebra $\fracg$ in $\Omega$. It follows from the universal defining property of $U(\fracg)$ that $L$ extends uniquely as a representation $h\mapsto L_h$ of the unital associative algebra $U(\fracg)$ in $\Omega$, i.e. as a homomorphism of $U(\fracg)$ into the algebra $\End^0(\Omega)$ of endomorphisms of degree 0 of $\Omega$. This extension will be refered to as {\sl the canonical extension of $L$ to $U(\fracg)$}.  

Let us recall that $U(\fracg)$ is not only an algebra but that it is a Hopf algebra with unique coproduct $\Delta$, counit $\varepsilon$ and antipode $S$ such that
\begin{equation}
\Delta X = X\otimes \bbbone +\bbbone \otimes X,\ \ \varepsilon(X)=0,\ \ S(X)=-X \qquad \text{for} \quad X\in \fracg .
\label{HU}
\end{equation}
 
\begin{proposition}\label{LextU}
The canonical extension of $L$ to $U(\fracg)$ has the following properties :\\
$\mathrm{(a)}$ $L_hd=dL_h$ for any $h\in U(\fracg)$\\
$\mathrm{(b)}$ $L_h(\bbbone)=\varepsilon(h)\bbbone$ for any $h\in U(\fracg)$ where $\bbbone$ is the unit of $\Omega$,\\
$\mathrm{(c)}$ $L_h(\alpha\beta)=\sum_i L_{h^{(1)}_i}(\alpha) L_{h^{(2)}_i}(\beta)$  
or any $h\in U(\fracg)$ 
and for any $\alpha,\beta \in \Omega$.\\
$\mathrm{(d)}$ $i_XL_h=\sum\nolimits_j L_{h_j^{(1)}} i_{\text{ad}(h^{(2)}_j)X}$ for any $X\in \fracg$ and $h\in  U(\fracg)$, where
the right adjoint action $\ad$ is defined by $
\ad(h)g=\sum\nolimits_j S(h_j^{(1)})g h_j^{(2)}$
for any $g,h\in U(\fracg)$
\end{proposition}

\noindent \underbar{Proof}. 
(a) is clear since $L_{\bbbone}$ is the identity mapping of $\Omega$ onto itself, the $L_X$ commute with $d$ for $X\in \fracg$ and since the unit $\bbbone\in  U(\fracg)$ and the $X\in \fracg$ generate $U(\fracg)$.  

\noindent \phantom{\underbar{Proof}.} (b) follows from the fact that $L_X (\bbbone)=0$ for $X\in \fracg$ so $L_h (\bbbone)=0$ for $h\in \ker(\varepsilon)$ and from $L_{\bbbone}(\bbbone)=\bbbone$ (since $L_{\bbbone}$ is the identity mapping of $\Omega$).  

\noindent \phantom{\underbar{Proof}.} (c) follows from the fact that it holds for $h=\bbbone (\in U(\fracg))$ and for $h=X\in \fracg$ and that furthermore both $h\mapsto L_h$ and $h\mapsto \Delta h$ are multiplicative homomorphisms.

\noindent \phantom{\underbar{Proof}.} (d) Relation (\ref{axop}) extends as the following identity
\begin{equation}
i_X L_{Y_1\dots Y_n}=\sum_p \sum_{(i,j)\in (p,n-p)\ \text{shuffles}} L_{Y_{i_1}\dots Y_{i_p}} i_{[[\dots[X, Y_{j_1}],\dots], Y_{j_{n-p}}]}
\label{Lig}
\end{equation}
for $X$ and the $Y_k$ in $\fracg$ which implies (d) by using the fact that one has 
\begin{equation}
\Delta (Y_1,\dots Y_n)=\sum_p \sum_{(i,j)\in (p,n-p)\text{shuffles}} Y_{i_1}\dots Y_{i_p} \otimes Y_{j_1} \dots Y_{j_{n-p}}
\label{Cop}
\end{equation}
for the coproduct of $Y_1\cdots Y_n$.~$\qed$\\

In contrast to $L$, $i$ has no canonical extension to $U(\fracg)$ which could be used to define more generally the axioms for an ``operation" of an arbitrary Hopf algebra $\calh$. Nevertheless, it is worth noticing here that $X\mapsto i_X$ admits an extension $h\mapsto \bar \imath_h$ to $U(\fracg)$ which satisfies
\begin{equation}
\bar\imath_{\bbbone}=0,\>\> \>\> L_h=d\bar\imath_h+\bar\imath_hd+\varepsilon(h)I_\Omega,\>\>\>\>\bar\imath_g L_h=\sum_jL_{h^{(1)}_j}\bar\imath_{\ad (h^{(2)}_j)g}
\label{propbari}
\end{equation}
for any $h,g\in U(\fracg)$. This extension is defined by setting
\begin{equation}
\bar\imath_{X^{n+1}}=L_{X^n}i_X=(L_X)^n i_X
\label{bari}
\end{equation}
for any $X\in\fracg$, $n\in \mathbb N$ and by using the Poincaré-Birkhoff-Witt theorem together with $\bar\imath_{\bbbone}=0$. However this extension $\bar\imath$, as well as the complicated rule generalizing the antiderivation property which it satisfies, depends in a crucial way on the particular structure of enveloping algebra, i.e. generation by primitive elements and the Poincaré-Birkhoff-Witt theorem. Let us just observe that this example shows that the properties (\ref{propbari}), together with the properties (a), (b), (c) of Proposition \ref{LextU}, which only depend on the general Hopf algebra structure are consistent (and natural).

\subsection{Algebraic connections in commutative $\fracg$-operations}\label{acgp}

We now give a short review of the notion introduced in \cite{car:1951a} of algebraic connection for operations of $\fracg$ in graded \underbar{commutative} differential algebras. These $\fracg$-operations will be refered to as {\sl commutative $\fracg$-operations}.

Given such an operation of a Lie algebra $\fracg$ in a graded \underbar{commutative} differential algebra $\Omega$ {\sl an algebraic connection} or simply  {\sl a connection} in $\Omega$ is a linear mapping 
\begin{equation}
\alpha:\fracg^\ast\rightarrow \Omega^1
\label{Ccl}
\end{equation}
of the dual vector space $\fracg^\ast$ of $\fracg$ such that, 
for any $X\in \fracg$ and $\theta\in \fracg^\ast$, one has
\begin{equation}
i_X(\alpha(\theta))=\theta(X)
\qquad \text{and} \qquad
L_X(\alpha(\theta))=\alpha(\theta\circ\ad(X)) \, .
\label{CclL}
\end{equation}

By the universal property of the exterior algebra, the mapping $\alpha$ extends as an homomorphism again denoted by 
\begin{equation}
\alpha:\wedge\fracg^\ast\rightarrow \Omega
\label{Cclex}
\end{equation}
of graded commutative algebras. In fact, $\wedge\fracg^\ast$ is a graded differential algebra (endowed with the Koszul differential $d$) and the {\sl curvature of} $\alpha$ is the linear mapping
\begin{equation}
\varphi:\fracg^\ast\rightarrow \Omega^2 \, , \qquad
\varphi(\theta)=(d\alpha-\alpha d)(\theta)
\label{Ccurv}
\end{equation}
for any $\theta\in \fracg^\ast$. Thus the curvature is the obstruction for $\alpha$ to be a homomorphism of graded differential algebras. It follows from the definitions that 
for any $X\in \fracg$ and $\theta\in \fracg^\ast$, one has
\begin{equation}
i_X(\varphi(\theta))=0
\qquad \text{and} \qquad
L_X(\varphi(\theta))=\varphi(\theta\circ \ad(X)) \, .
\label{InvC}
\end{equation}

In the standard example where $\Omega$ is the graded differential algebra of differential forms on a principal bundle $P(M,G)$ with  structure group $G$ such that $\mathrm{Lie}(G)=\fracg$, an algebraic connection on $\Omega$ is an ordinary principal bundle connection on $P(M,G)$.  

In Appendix B we give other standard notations for connections.

Let $\fracg$ be a fixed Lie algebra and consider the operations of $\fracg$ with connections in graded commutative differential algebras. There is a straightforward notion of morphism for such objects and one gets the category of {\sl commutative $\fracg$-operations with connections}. A morphism is a homomorphism of graded differential algebras which intertwins the $\fracg$-operations and which maps the connection on the connection.
There is a universal initial object $W(\fracg)$ in this category which is called the {\sl Weil algebra of the Lie algebra} $\fracg$. 
As graded commutative algebra 
\[
W(\fracg) = \wedge\fracg^\ast \otimes S\fracg^\ast
\]
is the tensor product of the exterior algebra 
$\wedge\fracg^\ast$ of the dual vector space $\fracg^\ast$ of $\fracg$ with the symmetric algebra $S\fracg^\ast$ where the degree $2n$ is given to the elements of $S^n\fracg^\ast$.  

Let $\theta$ be an element of $\fracg^\ast$, we denote by $\alpha(\theta)$ the element $\theta\otimes \bbbone$ of $W^1(\fracg)$ and by $\varphi(\theta)$ the element $\bbbone\otimes \theta$ of $W^2(\fracg)$. It is clear that there is a unique differential on $W(\fracg)$ for which 
\[
d\alpha(\theta)=\alpha(d\theta)+\varphi(\theta)
\]
for any $\theta\in \fracg^\ast$. In fact $W(\fracg)$ can be defined as well by a change of generators as the free graded commutative differential algebra generated by $\alpha(\fracg^\ast)$ in degree 1 and by $d\alpha(\fracg^\ast)$ in degree 2, which is a contractible algebra \cite{sul:1977}.  

The operation $X\mapsto i_X$ of $\fracg$ in $W(\fracg)$ is defined to be the unique antiderivation of $W(\fracg)$ such that 
$i_X(\alpha(\theta))=\theta(X)$
and 
$
i_X(\varphi(\theta))=0
$
for any $\theta\in \fracg^\ast$, ($X\in\fracg$). One verifies that all the axioms for an operation of $\fracg$ in the graded differential algebra $W(\fracg)$ are satisfied and that $\alpha$ is then a connection with curvature $\varphi$ in $W(\fracg)$.  

By the very definition of $W(\fracg)$, its cohomology is trivial and one can show by introducing the appropriate contracting homotopy that its invariant cohomology is also trivial. The basic differential subalgebra $W_B(\fracg)$ of $W(\fracg)$ is given by
$W^{2n}_B(\fracg)=\bbbone \otimes\cali^n(\fracg)$ and $W^{2n+1}_B(\fracg)=0$, 
where $\cali^n(\fracg)$ is the vector space of all ad-invariant homogeneous polynomials of degree $n$ on $\fracg$. It follows that $W_B(\fracg)$ coincides with its cohomology, i.e. the basic cohomology $H_B(W(\fracg))$.

The definition of $W(\fracg)$ implies that, given an operation of $\fracg$ with connection in a graded commutative differential algebra $\Omega$, there is a unique homomorphism of graded differential algebras of $W(\fracg)$ into $\Omega$ which is a morphism of commutative $\fracg$-operation with connection. It can be shown that this homomorphism induces in basic cohomology an homomorphism which does not depend on the connection of $\Omega$ but only depends on the operation of $\fracg$ in $\Omega$. This is the algebraic version of the Weil homomorphism. One recovers the familiar version by applying it to the standard example, remembering that the basic cohomology of the Weil algebra $W(\fracg)$ is isomorphic to the algebra $\cali(\fracg)$ of all ad-invariant polynomials on $\fracg$. Thus if $P(M,G)$ is a principal bundle over $M$ (the basis) with structure group $G$ such that $\fracg= \mathrm{Lie}(G)$, the Weil homomorphism is an algebra homomorphism from $\cali(\fracg)$ into the de Rham cohomology $H(M)$ of (the basis) $M$ such that the image of $\cali^n(\fracg)$ is contained in $H^{2n}(M)$ for any $n\in \mathbb N$.  

Let $P\in \cali(\fracg)$ be an $\ad$-invariant polynomial on $\fracg$. Then $\bbbone\otimes P\in W(\fracg)$ is closed and invariant so in view of the triviality of the invariant cohomology of $W(\fracg)$, one has $\bbbone\otimes P=dQ$ with $Q\in W_I(\fracg)$. Let 
$
\rho:W(\fracg)\rightarrow \wedge\fracg^\ast
$
 be the canonical projection. The image $\rho(Q)$ is an invariant form in $\wedge\fracg^\ast$ and it is not hard to show that it is independent of the choice of $Q$ as above. The corresponding linear mapping 
$
 \gamma:\cali(\fracg)\rightarrow \wedge_I\fracg^\ast
$
  from $\cali(\fracg)$ into the space of invariant forms on $\fracg$ is {\sl the Cartan map}. One has 
$
  \gamma(\cali^n(\fracg))\subset \wedge^{2n-1}_I \fracg^\ast
  $
  for any $n\geq 1$.

\section{Operations of Hopf algebras}

In this section we define the notion of operation of a Hopf algebra in a graded differential algebra.

\subsection{Definition of $\calh$-operations }\label{def}

Let $\calh$ be a Hopf algebra with coproduct $\Delta$, counit $\varepsilon$ and antipode $S$ and let $\Omega$ be a graded differential algebra with differential $d$.  
 {\sl An operation of the Hopf algebra $\calh$ in the graded differential algebra} $\Omega$ is a linear mapping $h\mapsto i_h$ of $\calh$ into the vector space $\End^{-1}(\Omega)$ of homogeneous linear endomorphisms of degree -1 of $\Omega$ satisfying 
\begin{equation}
i_{\bbbone}=0
\label{nor}
\end{equation}
and such that, by setting for $h\in \calh$
\begin{equation}
L_h=di_h+i_hd+\varepsilon (h) I_\Omega
\label{defLie}
\end{equation}
where $I_\Omega$ is the identity mapping of $\Omega$,  
the following properties hold. Firstly, for any $h\in \calh$ \begin{equation}
i_h(\alpha\beta)=\sum\nolimits_j i_{h^{(1)}_j}(\alpha)L_{h_j^{(2)}}(\beta)+(-1)^a\alpha \, i_h(\beta)
\label{antid}
\end{equation}
for $\alpha\in \Omega^a, \beta\in \Omega$. Secondly, for $g\in \calh$, 
\begin{equation}
i_g L_h=\sum\nolimits_j L_{h^{(1)}_j} i_{\text{ad}(h^{(2)}_j)g} 
\label{Cr}
\end{equation}
where the right adjoint action ad is defined (as before) by
\begin{equation}
\text{ad}(h)g=\sum\nolimits_j S(h^{(1)}_j)g h^{(2)}_j \qquad \forall h, g \in \calh 
\label{adR}
\end{equation}
and finally
\begin{equation}
L_h L_g=L_{hg} \qquad \forall h,g \in \calh \, .
\label{alghom}
\end{equation}

By using the associativity of $\Delta$, the axioms for the counit $\varepsilon$ and of the antipode $S$, one verifies that Equation (\ref{Cr}) is equivalent to
\begin{equation}
\sum\nolimits_j L_{S(h^{(1)}_j)} i_g  L_{h^{(2)}_j}=i_{\text{ad}(h)g}
\label{equiI}
\end{equation}
which implies, by using (\ref{defLie}) that
\begin{equation}
\sum\nolimits_j L_{S(h^{(1)}_j)} L_g  L_{h^{(2)}_j}=L_{\text{ad}(h)g}
\label{equiL}
\end{equation}
for any $g, h\in \calh$.   One sees that Equation (\ref{equiL}) is also implied by the axiom (\ref{alghom}).  
 Equations (\ref{equiI}) and (\ref{equiL}) mean the ad-equivariance of the mappings $h\mapsto i_h$ and $h\mapsto L_h$. From (\ref{defLie}) and (\ref{nor}) it follows that one has
\begin{equation}
L_{\bbbone}=I_\Omega
\label{gLun}
\end{equation}
and
\begin{equation}
L_h(\bbbone)=\varepsilon(h)\bbbone
\label{gL1}
\end{equation}
for  any $h\in \calh$, where in (\ref{gLun}) $\bbbone$ is the unit of $\calh$ while in (\ref{gL1}) $\bbbone$ is the unit of $\Omega$. From (\ref{defLie}) and (\ref{antid}) it follows that one has
\begin{equation}
L_h(\alpha\beta)=\sum\nolimits_j L_{h^{(1)}_j}(\alpha) L_{h_j^{(2)}}(\beta)
\label{cohom}
\end{equation}
for any $h\in \calh$ and $\alpha, \beta\in \Omega$.  
Finally, from (\ref{defLie}) it follows that one has
\begin{equation}
L_hd=dL_h
\label{cLd} \, 
\end{equation}
for any $h\in \calh$.

Operations of $\calh$ in graded differential algebras will be also refered to as $\calh$-{\sl operations}.
Given an operation $i$ of $\calh$ in $\Omega$ and an operation $i'$ of $\calh$ in $\Omega'$, a {\sl morphism of $\calh$-operation} from $(i,\Omega)$ to $(i',\Omega')$ is a homomorphism $f:\Omega\rightarrow \Omega'$ of graded differential algebra which satisfies
$f(i_h(\omega))=i'_h(f(\omega))$
for any $h\in \calh$ and $\omega\in \Omega$.
The $\calh$-operations and their morphisms form a category which will be refered to as {\sl the category of $\calh$-operations}.

\subsection{Invariance, horizontality and basicity}

Given an operation of $\calh$ in $\Omega$ as in \S\ref{def}, an element $\alpha \in \Omega$ will be said to be {\sl invariant} if 
\begin{equation}
L_h(\alpha)=\varepsilon(h)\alpha
\label{ginv}
\end{equation}
for any $h\in \calh$; $\alpha$ will be said to be {\sl horizontal} if 
\begin{equation}
i_h(\alpha)=0
\label{ghor}
\end{equation}
for any $h\in \calh$, and finally $\alpha$ will be said to be {\sl basic} if it is both invariant and horizontal.  

It follows from (\ref{defLie}) and (\ref{cohom}) that the set $\Omega_I$ of all invariant elements of $\Omega$ is a graded differential subalgebra of $\Omega$, it follows from (\ref{antid}) and (\ref{Cr}) that the set $\Omega_H$ of all horizontal elements of $\Omega$ is a graded subalgebra of $\Omega$ which is stable by the $L_h$ for $h\in \calh$ and it follows from (\ref{defLie}) again that the set $\Omega_B=\Omega_I\cap \Omega_H$ of all basic elements is a graded differential subalgebra of $\Omega_I$ and therefore also of $\Omega$. The cohomology $H_I(\Omega)$ of $\Omega_I$ will be refered to as the invariant cohomology of $\Omega$ while the cohomology $H_B(\Omega)$ of $\Omega_B$ will be refered to as the basic cohomology of $\Omega$ for the operation of $\calh$ in $\Omega$.

\subsection{Operations of the Lie algebra $\fracg$ and of the Hopf algebra $U(\fracg)$}

Let $\fracg$ be a Lie algebra and let $h\mapsto i_h$ be an operation of the Hopf algebra $U(\fracg)$ in the graded differential algebra $\Omega$. It is clear that by restriction to $\fracg\subset U(\fracg)$ one obtains an operation $X\mapsto i_X$ of the Lie algebra $\fracg$ in the graded differential algebra $\Omega$, (in the sense of \S \ref{defope}). For the converse, one has the following result.

\begin{proposition}\label{gU}
Let $\Omega$ be a differential envelope, that is $\Omega=\Omega(\cala)$ for some algebra $\cala$. Then an operation of the Lie algebra $\fracg$ in $\Omega$ has a unique extension as an operation of the Hopf algebra $U(\fracg)$ in $\Omega$.
\end{proposition}

\noindent \underbar{Proof}. Any element of $\Omega$ is a linear combination of terms of the form
\[
x_0dx_1 \dots dx_n
\]
with $x_\alpha \in \Omega^0$. Let $h\mapsto i_h$ be an operation of $U(\fracg)$ in $\Omega$. One has
\[
i_h(x_0dx_1\dots dx_n)=x_0 i_h(dx_1\dots dx_n)
\]
and $i_h(dx_\alpha)=L_h(x_\alpha)-\varepsilon (h) x_\alpha$. It follows that the operation is completely specified by the $L_h$ in view of (\ref{antid}). On the other hand $L_h$ is completely specified by the $L_X$ for $X\in \fracg$ in view of (\ref{alghom}) since $\fracg$ generates $U(\fracg)$. Thus, starting from the $i_X$ with $X\in \fracg$ one constructs the $L_h$ for $h\in U(\fracg)$ and the $i_h$ for $h\in U(\fracg)$ by using the above formulae. One verifies that the axioms of \S \ref{def} are satisfied. $\qed$

If $\Omega$ is generated in degree 0 as graded differential algebra, the above proof shows that if an operation of $\fracg$ in $\Omega$ has an extension as an operation of $U(\fracg)$ in $\Omega$, then this extension is unique. Furthermore under these conditions the notions of invariance, horizontality and basicity are the same for the operation of the Lie algebra $\fracg$ and for the operation of the Hopf algebra $U(\fracg)$. However, it is worth noticing here that there exist $\fracg$-operations which do not admit extensions as $U(\fracg)$-operations. For instance Property (\ref{antid}) for $U(\fracg)$-operations is not compatible with the (eventual) graded commutativity of $\Omega$, i.e. $U(\fracg)$-operations cannot be commutative. Indeed, by applying (\ref{antid}), one gets
\[
i_h(fd(g)-d(g)f)=fL_h(g)+gL_h(f)-\varepsilon(h)fg-L_h(gf)
\]
for $f,g\in \Omega^0$, which implies
\[
i_{X^2}(fd(g)-d(g)f)=-2L_X(g)L_X(f)+[f,L_{X^2}(g)]
\]
for $X\in \fracg$, which is incompatible with the commutativity whenever $L_X(g)L_X(f)\not=0$. Notice also that $\Omega=\wedge\fracg^\ast$ and, more generally, the graded differential algebras associated by Koszul duality to Lie prealgebras \cite{mdv-lan:2011} are not generated in degree 0 as graded differential algebras.

\subsection{Filtration of a $\calh$-operation}

Assume that one has an operation of $\calh$ in $\Omega$. Then there is an associated decreasing filtration $F^p(\Omega)$ of $\Omega$ given by setting
\begin{equation}
F^p(\Omega^n)=\{\omega\in \Omega^n\vert i_{h_1}\cdots i_{h_{n-p+1}}(\omega)=0,\>\> \forall h_k\in\calh\}
\label{filtr}
\end{equation}
for $n\geq p$ and $F^p(\Omega^n)=0$ for $n<p$. One has
\begin{equation}
i_h(F^p(\Omega))\subset F^p(\Omega),\>\> \forall h\in \calh
\label{sti}
\end{equation}
\begin{equation}
d(F^p(\Omega))\subset F^p(\Omega)
\label{std}
\end{equation}
which implies
\begin{equation}
L_h(F^p(\Omega)) \subset F^p(\Omega)
\label{stL}
\end{equation}
for any $p\in \mathbb N$, and
\begin{equation}
F^p(\Omega)F^q(\Omega)\subset F^{p+q}(\Omega)
\label{coher}
\end{equation}
for  any $p,q\in \mathbb N$ (with $F^{p+1}(\Omega)\subset F^p(\Omega)$ and $F^0(\Omega)=\Omega$).
To such a filtration of graded differential algebra of $\Omega$ corresponds a convergent spectral sequence $(E_r(\Omega),d_r)_{r\in \mathbb N}$ where $E_r(\Omega)$ is a bigraded algebra $E_r(\Omega)=\oplus_{p,q\in \mathbb N} E^{p,q}_r(\Omega)$ and where $d_r$ is an homogeneous differential on $E_r(\Omega)$ of bidegree $(r,1-r)$.
It is clear that $\omega\in F^p(\Omega^{p+q})$ is equivalent to $i_{h_1}\cdots i_{h_q}(\omega)\in \Omega^p_H$ for any $h_k\in \calh$ and therefore $\omega\mapsto i_{h_1}\cdots i_{h_q}(\omega)$ defines a linear mapping $\varphi^{p,q}$ of $F^p(\Omega^{p+q})$ into the linear space $C^q(\calh,\Omega^p_H)$ of Hochschild $q$-cochains of $\calh$ with values in $\Omega^p_H$ equipped with the left $\calh$-action $h\mapsto L_h$ and the trivial right action given by the counit $\varepsilon$. The kernel of $\varphi^{p,q}$ is $F^{p+1}(\Omega^{p+q})$ so one has the exact sequences
\[
0\rightarrow F^{p+1}(\Omega^{p+q})\stackrel{\subset}{\rightarrow}F^p(\Omega^{p+q})\stackrel{\varphi^{p,q}}{\rightarrow}C^q(\calh,\Omega^p_H)
\]
of vector spaces for $p,q\in \mathbb N$. By setting
\[
E^{p,q}_0(\Omega)=F^p(\Omega^{p+q})/F^{p+1}(\Omega^{p+q})
\]
one has $E^{p,q}_0(\Omega)\subset C^q(\calh,\Omega^p_H)$ and thus the 0-term of the spectral sequence associated with the filtration
$E_0(\Omega)=\oplus_{p,q} E^{p,q}_0(\Omega)$
is a bigraded subalgebra of $\oplus_{p,q} C^q(\calh,\Omega^p_H)$. Concerning the first term $E_1(\Omega)=H(E_0(\Omega))$ of the spectral sequence, one has
$E^{p,0}_0(\Omega)=F^p(\Omega^p)$ and therefore 
$E^{p,0}_1(\Omega)=H^{p,0}(E_0(\Omega))=\{\omega\in F^p(\Omega^p)\vert d\omega\in F^{p+1}(\Omega^{p+1})\}$
for any $p\in \mathbb N$. Thus one has 
$E^{p,0}_1(\Omega)=\{\omega\in\Omega^p\vert i_h(\omega)=0\>\> \text{and}\>\>i_h d (\omega)=0,\>\> \forall h\in\calh\}$
which is equivalent to
\begin{equation}
E^{p,0}_1(\Omega)=\Omega^p_B
\label{sbasis}
\end{equation}
for $p\in \mathbb N$, so $E^{p,0}_2(\Omega)=H^p_B(\Omega)$, etc.

\section{Theory of connections}

In this section we introduce and study a noncommutative generalization of the theory of algebraic connections on the operations of a Lie algebra \cite{car:1951b}, \cite{gre-hal-van:1976}.

\subsection{The graded differential algebra $C(\calh)$}

Although in the sequel $\calh$ is a Hopf algebra, in this subsection only its algebra structure is involved. At the end of it the augmentation $\varepsilon$ (counits) of $\calh$ will also play a role.  

Let $C(\calh)=\oplus_n C^n(\calh)$ be the graded algebra of multilinear forms on $\calh$, that is one has
\[
C^n(\calh)=(\calh^{\otimes^n})^\ast
\]
and the product is the tensor product of multilinear forms.
The product $\mu:\calh\otimes \calh\rightarrow \calh$ of $\calh$ induces by transposition the linear mapping
\[
\mu^t:\calh^\ast\rightarrow (\calh\otimes \calh)^\ast
\]
so $-\mu^t$ (the minus sign is here to match the usual convention) is a linear mapping of $C^1(\calh)$ into $C^2(\calh)$ which has an extension
\[
d_0:C(\calh)\rightarrow C(\calh)
\]
as an antiderivation of degree 1 of $C(\calh)$ given by
\begin{equation}
d_0(\Psi)(h_0,h_1,\dots,h_n)=\sum_{1\leq k\leq n} (-1)^k \Psi(h_0,\dots,h_{k-2},h_{k-1}h_k,h_{k+1},\dots,h_n)
\label{dC}
\end{equation}
for $\Psi\in C^n(\calh)$ and $h_0,\dots,h_n\in \calh$. The associativity of the product of $\calh$ is equivalent to
\begin{equation}
d_0^2=0
\label{d2C}
\end{equation}
so that ($C(\calh),d_0$) is a graded differential algebra. The construction of the graded differential algebra $(C(\calh),d_0)$ works as well for any unital associative algebra, it is dual to the acyclic bar complex, and it is a standard fact that its cohomology vanishes in positive degrees.  

The counit $\varepsilon$ gives a structure of $\calh$-bimodule to the ground field $\mathbb K$, 
referred to as {\sl the trivial bimodule} $\mathbb K$. Thus $C(\calh)$ is also the space of Hochschild cochains of $\calh$ with coefficients in the trivial bimodule $\mathbb K$. The corresponding Hochschild differential $d$ reads
\begin{equation}
d\omega=d_0\omega+\varepsilon\omega + (-1)^{n+1}\omega\varepsilon  \qquad \text{for} \quad \omega\in C^n(\calh) \, . 
\label{Hoch}
\end{equation}
In contrast to the cohomology of $C(\calh)$ for $d_0$, the cohomology of $d$ is nontrivial in general. For instance, in the case $\calh=U(\fracg)$, this cohomology is isomorphic to
the cohomology of the Lie algebra $\fracg$, see e.g. in \cite{lod:1992}. In fact one defines a quasi-isomorphism
\begin{equation}
C(U(\fracg))\rightarrow \wedge\fracg^\ast
\label{quasiU}
\end{equation}
of $(C(U(\fracg)),d)$ onto $\wedge\fracg^\ast$ endowed with the Koszul differential, by restriction to $\fracg$ (or more precisely to $T(\fracg)$) and antisymmetrization. Thus, from this point of view $C(\calh)$ equipped with the differential $d$ is the analogue of $\wedge\fracg^\ast$ equipped with the Koszul differential. In the following $C(\calh)$ endowed with the differential $d$ will be refered to as {\sl the graded differential algebra} $C(\calh)$.

\subsection{The operation of $\calh$ in $C(\calh)$} \label{ohinch}

Motivated by the requirement \eqref{antid}, one defines an operation of the Hopf algebra $\calh$ in the graded differential algebra $C(\calh)$ in the following manner. 
Let $\Psi\in C^n(\calh)$ be a $n$-linear form on $\calh$ and let us define, for $h\in \calh$, $i_h(\Psi)\in C^{n-1}(\calh)$  by
\begin{equation}
\begin{array}{l}
 i_h(\Psi)(g_1,\dots,g_{n-1}) =
 \\
 \\
\displaystyle{\sum^{n-2}_{p=0}}(-1)^p \sum_{i_p}\Psi(g_1,\dots,g_p,h^{(1)}_{i_p}-\varepsilon (h^{(1)}_{i_p})\bbbone, \text{ad}(h^{(2)}_{i_p})g_{p+1},\dots,\text{ad}(h_{i_p}^{(n-p)})g_{n-1})\\
\\
+ (-1)^{n-1}\Psi(g_1,\dots,g_{n-1},h-\varepsilon(h)\bbbone)
\end{array}
\label{intC}
\end{equation}
for $g_k\in \calh$. Here we have set
\begin{equation}
(\Delta\otimes I_\calh^{\otimes^{n-p-2}})\dots (\Delta\otimes I_\calh)\Delta h=\sum\nolimits_{i_p}h^{(1)}_{i_p}\otimes h^{(2)}_{i_p}\otimes \dots \otimes h^{(n-p)}_{i_p}
\label{Itcop}
\end{equation}
for the iterated coproducts (which occur of course only for $n\geq 2$).  

One verifies that this defines a $\calh$-operation and that $L_h=di_h+i_hd+\varepsilon_h$ is given by

\begin{equation}
L_h(\Psi)(g_1,\dots,g_n)=\sum\nolimits_{i_0}\Psi(\text{ad}(h^{(1)}_{i_0})g_1,\dots,\text{ad}(h^{(n)}_{i_0})g_n)
\label{LieC}
\end{equation}
with obvious notations. One also verifies that
\begin{equation}
i_hd+di_h=i_hd_0+d_0i_h
\label{dd0}
\end{equation}
which implies that $h\mapsto i_h$ is an operation of $\calh$ in the graded differential algebra $(C(\calh),d_0)$ as well. 
The invariant cohomology of $(C(\calh),d_0)$  also vanishes in positive degrees.

\subsection{Algebraic connections}\label{alco}

Let $\calh$ be a Hopf algebra, $\Omega$ be a graded differential and assume that one has an operation $h\mapsto i_h$ of $\calh$ in $\Omega$, (i.e. that 
$(i,\Omega)$ is a $\calh$-operation).  

{\sl An algebraic connection on the $\calh$-operation $(i,\Omega)$} or simply {\sl a connection on} $(i,\Omega)$ is a homomorphism of graded algebras 
\begin{equation}
\alpha:C(\calh)\rightarrow \Omega
\label{dCo0}
\end{equation}
such that, for any $\Psi\in C(\calh)$, it satisfy 
\begin{equation}
i_h(\alpha(\Psi))=\alpha(i_h(\Psi))
\qquad \text{and} \qquad
L_h(\alpha(\Psi))=\alpha(L_h(\Psi)) \, .
\label{GCo2}
\end{equation}
In particular, for $\psi\in \calh^\ast$ this implies that, for any $h\in \calh$, 
\begin{equation}
i_h(\alpha(\psi))=\psi(h)-\varepsilon(h)\psi(\bbbone) \qquad \text{and} \qquad 
L_h(\alpha(\psi))=\alpha(\psi\circ \text{ad}(h)) \, .
\label{dCo2}
\end{equation}

The {\sl curvature of $\alpha$} is the homogeneous linear mapping of degree 1
\[
\varphi:C(\calh)\rightarrow \Omega
\]
defined by
\begin{equation}\label{courb}
\varphi(\Psi) =(d\alpha-\alpha d)(\Psi)
\end{equation}
for any $\Psi\in C(\calh)$. As in the classical theory, the curvature is the obstruction for $\alpha$ to be a homomorphism of graded differential algebras. This implies that
\begin{equation}
i_h(\varphi(\psi))=0   \qquad \text{and} \qquad 
L_h(\varphi(\psi))=\varphi(\psi\circ \text{ad}(h))
\label{dequiF}
\end{equation}
for any $h\in \calh$ and $\psi\in \calh^\ast$, and more generally
\begin{equation}
i_h(\varphi(\Psi))=-\varphi(i_h(\Psi))
\qquad \text{and} \qquad
L_h(\varphi(\Psi))=\varphi(L_h(\Psi))
\label{GF2}
\end{equation}
for any $\Psi\in C(\calh)$, together with 
\begin{equation}
\varphi(\bbbone)=0
\label{normF}
\end{equation}
that is $\varphi(C^0(\calh))=0$. Furthermore (\ref{courb}) implies also that
\begin{equation}
\varphi(\Phi\Psi)=\varphi(\Phi)\alpha(\Psi)+(-1)^f\alpha(\Phi)\varphi (\Psi)
\label{Fder}
\end{equation}
for $\Phi\in C^f(\calh), \Psi\in C(\calh)$ and that one has the following version of Bianchi identity
\begin{equation}
d(\varphi(\Psi))=-\varphi(d\Psi) \qquad \text{for} \quad \Psi\in C(\calh) \, .
\label{Bianchi}
\end{equation} 

Notice that $\alpha$ induces a structure of $C(\calh)$-bimodule on $\Omega$ and that then (\ref{Fder}) means that $\varphi$ is an antiderivation of $C(\calh)$ into $\Omega$, (see in \S \ref{dega}).  

One verifies that the identity mapping  $I_{C(\calh)}$ of $C(\calh)$ onto itself is a connection $\alpha_C$ on the $\calh$-operation $C(\calh)$ defined in the last section and that this connection is {\sl flat}, that is has a vanishing curvature $(\varphi_C=0)$. This connection will be refered to as {\sl the canonical flat connection of $C(\calh)$}. It is worth noticing here that in the case $\calh=U(\fracg)$ where $\fracg$ is the Lie algebra of a Lie group $G$, $\fracg=\Lie (G)$, the Maurer-Cartan equation for the canonical invariant connection on $G$ refer directly to the flatness of this connection on $C(\calh)$.  

\noindent \underbar{Remark}:
With the Koszul convention the relations (\ref{GCo2})
and (\ref{GF2}) all appear as ``commutation properties", therefore by
using corresponding graded commutators, they read
$[\alpha, i_h]=0,\>\> [\alpha,L_h]=0$ and $[\varphi,i_h]=0,\>\> [\varphi, L_h]=0$
while (\ref{courb}) and (\ref{Bianchi}) can be written $\varphi=[d,\alpha]$ and $[d,\varphi]=0$.

\subsection{The set of connections on a $\calh$-operation}

 Let $\calh$ be a vector space and let $C(\calh)=\oplus_{n\geq 0} C^n(\calh)$ be the graded connected algebra of all multilinear forms on $\calh$. Let us consider the quotient
\[
V(\calh)=C^+(\calh)/(C^+(\calh))^2
\]
where $C^+(\calh)=\oplus_{n\geq 1} C^n(\calh)$. The product of $C^+(\calh)$ induces the trivial zero product on $V(\calh)$ which is now just a graded vector space. By choosing a 
supplementary of $(C^+(\calh))^2$ in $C^+(\calh)$, one obtains the following lemma.

\begin{lemma} \label{injlin}
There exists an injective homogeneous linear mapping of $V(\calh)$ into $C(\calh)$.
\end{lemma}

Let $\fracT(V(\calh))$ be the tensor algebra $T(V(\calh))$ of $V(\calh)$ endowed with the unique graduation of algebra which induces the graduation of $V(\calh)$. One has the following result.

\begin{proposition}\label{isoma}
The graded algebras $C(\calh)$ and $\fracT(V(\calh))$ are isomorphic.
\end{proposition}

\noindent\underline{Proof}. By Lemma \ref{injlin} there is an injective homogeneous linear mapping of $V(\calh)$ into $C(\calh)$. In view of the universal property of $\fracT(V(\calh))$ this mapping induces an homomorphism of the graded algebra $\fracT(V(\calh))$ into $C(\calh)$. It is easy to verify that this homomorphism is an isomorphism of graded vector spaces. $\qed$

\begin{proposition}\label{homlin}
The set of homomorphisms of graded algebras of $C(\calh)$ into a graded algebra $\Omega$ is in bijection with the set of homogeneous linear mappings of $V(\calh)$ into $\Omega$.
\end{proposition}

\noindent\underbar{Proof}. This is just the universal property of $\fracT(V(\calh))$ combined with Proposition \ref{isoma}.\\

This last proposition has the following obvious corollary.

\begin{corollary}\label{Affc}
Let $\calh$ be a Hopf algebra. Then the set of connections on a $\calh$-operation $\Omega$ admits a structure of affine space modeled on a vector subspace of the vector space of linear homomorphisms of graded vector spaces of $V(\calh)$ into $\Omega$.
\end{corollary}

The vector subspace in Corollary \ref{Affc} is the subspace of the graded linear homomorphisms
\[
a=\sum_{n\geq 1} a_n:V(\calh)\rightarrow \Omega
\]
such that the homogeneous part of the equations induced by (\ref{GCo2}) are satisfied. The inhomogeneous part which concerns only $a_1$ corresponds to the first equation of (\ref{dCo2}) which should be replaced by $i_h(a_1(\psi))=0$ for $\psi\in \calh^\ast$.\\

\noindent\underbar{Remark}. When $\calh$ is finite-dimensional one has $V(\calh)=\calh^\ast$, so the linear mapping of Lemma~\ref{injlin}, the isomorphism of Proposition \ref{isoma} and the bijection of Proposition \ref{homlin} are unique. This implies that the affine structure of Corollary \ref{Affc} is unique for a finite-dimensional Hopf algebra.

\section{The universal $\calh$-operation with connection $W(\calh)$}\label{uhowc}

In this section $\calh$ is a fixed Hopf algebra and we define a noncommutative version of the Weil algebra, the Weil algebra $W(\calh)$ of the Hopf algebra $\calh$.

\subsection{The category of $\calh$-operations with connections}

Let $f$ be a morphism of $\calh$-operation from $(i,\Omega)$ to $(i',\Omega')$ and let 
\[
\alpha:C(\calh)\rightarrow \Omega
\]
 be a connection on $(i,\Omega)$, then the image $f\circ \alpha$ of $\alpha$ by $f$
\[
f\circ \alpha : C(\calh)\rightarrow \Omega'
\]
is clearly a connection on $(i',\Omega')$ which will be denoted by $f(\alpha)$.  

Then, one defines the category of $\calh$-operations with connections in a natural way.
{\sl A $\calh$-operation with connection} $(i,\Omega,\alpha)$ is an $\calh$-operation $(i,\Omega)$ equipped with a connection $\alpha$. Given two $\calh$-operations with connections $(i,\Omega,\alpha)$ and $(i',\Omega',\alpha')$, {\sl a morphism of $\calh$-operation with connection} from $(i,\Omega,\alpha)$ to $(i',\Omega',\alpha')$ is a morphism of $\calh$-operation $f$ from $(i,\Omega)$ to $(i',\Omega')$ such that $\alpha'=f(\alpha)$.

 It turns out that this category of $\calh$-operations with connections has a universal initial object $W(\calh)$ which is the appropriate generalization of the Weil algebra in our context and which will be described in this section.
 We first describe, in the next subsection, the construction of the universal differential calculus over a graded algebra.
 
 \subsection{Differential envelopes of graded algebras}\label{dega}
 
 Let $\cala=\oplus_n\cala^n$ be a graded algebra ($\mathbb N$-graded, unital, associative) with product
 \[
 (x,y)\mapsto m(x\otimes y)=xy
 \]
 for $x,y\in \cala$. 
An {\sl antiderivation} of $\cala$ into a $\cala$-bimodule $\calm$ is a linear mapping
 \[
 \delta:\cala\rightarrow \calm
 \]
 such that one has
 \begin{equation}
 \delta(ab)=\delta(a)b+(-1)^r a \delta(b)
 \label{antider}
 \end{equation}
 for $a\in \cala^r$ and $b\in \cala$. Following \cite{coq-kas:1989} let us define the twisted product
$
 \mu:\cala\otimes \cala\rightarrow \cala
$
 on $\cala$ by setting
 \begin{equation}
\mu(a\otimes b)=(-1)^sab
\label{twpr}
\end{equation}
for $a\in \cala$ and $b\in \cala^s$. There is a structure of $\cala$-bimodule on $\cala\otimes \cala$ given by setting 
\[
x(a\otimes b)=xa\otimes b,\>\>\> (a\otimes b)y=a\otimes by
\]
for $x,y,a,b\in \cala$ and the kernel $J$ of $\mu$ is a sub-bimodule of $\cala\otimes\cala$. One verifies that one defines an antiderivation $d:\cala\rightarrow J$ of $\cala$ into the $\cala$-bimodule $J$ by setting
\begin{equation}
d(x)=\bbbone \otimes x - (-1)^n x\otimes \bbbone
\label{uantider}
\end{equation}
for $x\in \cala^n$. This antiderivation $d:\cala\rightarrow J$ is characterized (up to isomorphisms) by the following universal property.

\begin{proposition}\label{UanD}
Let $\delta:\cala\rightarrow\calm$ be an antiderivation of $\cala$ into an $\cala$-bimodule $\calm$. Then there is a unique homomorphism  of $\cala$-bimodules $i_\delta:J\rightarrow \calm$ of $J$ into $\calm$ such that $\delta=i_\delta\circ d$.
\end{proposition}

Thus this construction of \cite{coq-kas:1989} gives the counterpart for antiderivations of the classical construction \cite{car-eil:1973} of the universal derivations (see also \cite{bou:1970}), furthermore the proof of Proposition \ref{UanD} is the same as the proof of the corresponding proposition for derivations in \cite{car-eil:1973}, \cite{bou:1970}. A key remark for the proof is that as left $\cala$-module, one has the isomorphisms
\[
J\simeq \cala d\cala \simeq \cala\otimes d\cala\simeq \cala\otimes (\cala/\mathbb K \bbbone)
\]
while the right $\cala$-module structure is obtained from above by the graded Leibniz rule (\ref{antider}).  

The above apply as well for $\mathbb Z$-graded or $\mathbb Z_2$-graded algebras, in fact the paper \cite{coq-kas:1989} is written in the $\mathbb Z_2$-graded context. Nevertheless for the following $\cala$ is taken to be $\mathbb N$-graded.  

We now introduce a graduation on $J\subset \cala\otimes \cala$ by setting $J=\oplus_{n\geq 0}J^{n+1}$ with 
\[
J^{n+1}\subseteqq\oplus_{r+s=n}\cala^r\otimes \cala^s
\]
for $n\in \mathbb N$. Endowed with this graduation $J$ becomes a graded $\cala$-bimodule which will be denoted by $\Omega^1_{gr}(\cala)$. Thus
\[
d:\cala\rightarrow \Omega^1_{gr}(\cala)
\]
is a graded derivation of degree 1 of $\cala$ into $\Omega^1_{gr}(\cala)$. Note that the kernel of $d$ is $\mathbb K\bbbone$ so that
\[
d\cala\simeq (\cala/\mathbb K\bbbone)^{\bullet+1}
\]
i.e. $d\cala$ is $\cala/\mathbb K\bbbone$ with a shift +1 in graduation so $\Omega^1_{gr}(\cala)\simeq \cala\otimes (\cala/\mathbb K \bbbone)^{\bullet+1}$.  

Let us now define the graded algebra $\Omega_{gr}(\cala)=\oplus_n\Omega^n_{gr}(\cala)$ to be the tensor algebra over $\cala$ of the bimodule $\Omega^1_{gr}(\cala)$ endowed with the unique graduation of algebra which induces on $\cala=\Omega^0_{gr}(\cala)$ and on $\Omega^1_{gr}(\cala)$ their original graduation.  

The graded derivation $d$ of $\cala$ into $\Omega^1_{gr}(\cala)$ has a unique extension as a differential on $\Omega_{gr}(\cala)$, i.e. as a graded derivation of degree 1 of $\Omega_{gr}(\cala)$, again denoted by $d$, satisfying $d^2=0$. Endowed with this differential, $\Omega_{gr}(\cala)$ is a graded differential algebra which is characterized, up to an isomorphism, by the following universal property which is a graded counterpart of the universal property of the usual universal differential calculus $\Omega(\cala)$ over a non graded algebra $\cala$ \cite{ac:1982}, 
\cite{ac:1994}, \cite{kar:1983}, \cite{kar:1987}.

\begin{theorem}\label{UNIGR}
Any homomorphism of graded algebras
\[
\alpha:\cala\rightarrow \Omega
\]
of $\cala$ into a graded differential algebra $\Omega$ has a unique extension as homomorphism
\[
\Omega_{gr}(\alpha):\Omega_{gr}(\cala)\rightarrow \Omega
\]
of graded differential algebras.
\end{theorem}
The proof is completely similar to the proof of the ``ungraded counterpart".

Notice that by considering a non graded algebra $\cala$ as a graded algebra concentrated in degree 0 one can write $\Omega(\cala)=\Omega_\gr(\cala)$.
 
\subsection{The Weil algebra $W(\calh)$ of the Hopf algebra $\calh$}

Let $(i,\Omega,\alpha)$ be a $\calh$-operation with connection. Since then $\alpha$ is in particular a homomorphism of graded algebra of $C(\calh)$ into $\Omega$, it is natural to introduce the graded differential algebra
\begin{equation}
W(\calh)=\Omega_{gr}(C(\calh))
\label{grdifW}
\end{equation}
with the notations of the last subsection. With this, Theorem \ref{UNIGR} has the following corollary.

\begin{corollary}\label{UDifW}
Let $(i,\Omega,\alpha)$ be a $\calh$-operation with connection, then 
$\alpha:C(\calh)\rightarrow \Omega$
 has a unique extension
\[
W(\alpha):W(\calh)\rightarrow \Omega
\]
as homomorphism of graded differential algebras.
\end{corollary}

Let us denote by
\[
\alpha_W:C(\calh)\rightarrow W(\calh)
\]
the canonical injection of $C(\calh)$ into $W(\calh)$ as graded subalgebra and by
\[
\varphi_W=d\circ \alpha_W-\alpha_W\circ d: C(\calh)\rightarrow W(\calh)
\]
the corresponding obstruction for $\alpha_W$ to be a homomorphism of graded differential algebra. The homogeneous linear mapping $\varphi_W$ satisfies the relations (\ref{normF}), (\ref{Fder}) and (\ref{Bianchi}) that is
$\varphi_W(\bbbone)=0$ and, for $\Phi\in C^f(\calh)$, $\Psi\in C(\calh)$,  
\begin{align*}
\varphi_W(\Phi\Psi) &=\varphi_W(\Psi)\alpha_W(\Psi)+(-1)^f \alpha_W(\Phi)\varphi_W(\Psi) \\
~\\
d\varphi_W(\Psi) &=-\varphi_W(d\Psi) \, .
\end{align*}

\begin{proposition}\label{OpW}
There is a unique operation $h\mapsto i_h$ of the Hopf algebra $\calh$ in the graded differential algebra $W(\calh)$ for which $\alpha_W$ is an algebraic connection
\end{proposition}

\noindent\underbar{Proof}. In view of the definitions of \S \ref{alco}, if $\alpha_W$ is a connection, its curvature is given by $\varphi_W$. Thus one should have
$
i_h(\alpha_W(\Psi))=\alpha_W(i_h(\Psi))
$
and 
$
i_h(\varphi_W(\Psi))=-\varphi_W(i_h(\Psi))
$
for $h\in \calh$, $\Psi\in C(\calh)$.
This fixes the $i_h$ on $\alpha_W(C(\calh))\simeq C(\calh)$ and on $\varphi_W(C(\calh))$. One verifies that the relations $[\alpha_W,L_h]=0$ and $[\varphi_W,L_h]=0$ are satisfied and then the relation (\ref{antid}) fixes $i_h$ on $W(\calh)$. One then verifies that the so defined $h\mapsto i_h$ is an operation of $\calh$ in $W(\calh)$. $\qed$

By combining Corollary \ref{UDifW} and Proposition \ref{OpW} one arrives at the following theorem.
\begin{theorem}\label{inW}
Let $\Omega$ be a $\calh$-operation with connection, then there is a unique morphism of $\calh$-operation with connection from $W(\calh)$ to $\Omega$.
\end{theorem}

In other words $W(\calh)$ is a universal initial object in the category of $\calh$-operations with connections. As such it is unique up to isomorphism.  

The graded differential algebra $W(\calh)$ endowed with the structure described above will be refered to as  {\sl the Weil algebra of the Hopf algebra $\calh$.} 

It is clear that $W(\calh)$ plays in the present setting the same role as the Weil algebra $W(\fracg)$ of the Lie algebra $\fracg$ in the classical theory and it is worth reminding here that the role of $\wedge \fracg^\ast$ in the classical theory is played in our noncommutative framework by $C(\calh)$.  

One has canonically
\begin{equation}
W(\calh)=C(\calh)\oplus W_\varphi(\calh)
\label{decC}
\end{equation}
where $W_\varphi(\calh)$ is the two-sided ideal of $W(\calh)$ generated by $\varphi_W(C(\calh))$. This ideal is in fact a graded differential ideal of $W(\calh)$ so the corresponding canonical surjective homomorphism
\begin{equation}
\rho:W(\calh)\rightarrow C(\calh)
\label{canC}
\end{equation}
is a homomorphism of graded differential algebras. It is easy to see that it is the morphism of $\calh$-operation with connection of Theorem (\ref{inW}) for $\Omega=C(\calh)$ where the graded differential algebra $C(\calh)$ is endowed with its structure of $\calh$-operation with (flat) connection described in \S \ref{ohinch} and \S \ref{alco}. One has
\begin{equation}
\rho\circ \alpha_W=I_{C(\calh)}
\label{secC}
\end{equation}
which means that the connection of $W(\calh)$ is a section of $\rho$. Notice that for $\calh=U(\fracg)$, one obtains a surjective homomorphism
\begin{equation}
W(U(\fracg))\rightarrow W(\fracg)
\label{WUg}
\end{equation}
which is the counterpart of the quasi-isomorphism (\ref{quasiU}), by restriction of the arguments to $\fracg$ followed by graded symmetrization.

\subsection{Cohomology and invariant cohomology of $W(\calh)$}

The cohomology and the invariant cohomology of $W(\calh)$ are given by the following theorem.

\begin{theorem}\label{CoIW}
The cohomology $H(W(\calh))$ and the invariant cohomology $H_I(W(\calh))$ of $W(\calh)$ are both trivial, that is one has
\[
H^n(W(\calh))=H^n_I(W(\calh))=0
\]
for $n\geq 1$ while $H^0(W(\calh))$ and $H^0_I(W(\calh))$ identify to the ground field $\mathbb K$.
\end{theorem}

\noindent \underbar{Proof}. By its very definition, $W(\calh)$ is 
\[
W(\calh)=T_{C(\calh)}(\Omega^1_{gr}(C(\calh))) \, 
\]
the tensor algebra over $C(\calh)$ of the $C(\calh)$-bimodule $\Omega^1_{gr}(C(\calh))$.
Therefore, there is a unique antiderivation $K$ of $W(\calh)$ which is such that
\begin{equation}
K\circ \alpha_W=0
\label{Clin}
\end{equation}
and which satisfies
\begin{equation}
K\circ d\circ \alpha_W=\deg\circ \, \alpha_W \, ,
\label{KdH}
\end{equation}
where $\deg$ is the degree, (since $\alpha_W(C(\calh))$ and $d(\alpha_W(C(\calh)))$ generate $W(\calh)$).
In fact $\deg$ is a derivation of $W(\calh)$ into itself as well as a derivation of $C(\calh)$ into itself and one has 
$\deg \circ \, \alpha_W=\alpha_W\circ \deg$.  

 Notice that $\alpha_W(C(\calh)$ and $\varphi_W(C(\calh))$  generate as well $W(\calh)$ and (\ref{KdH}) is equivalent to 
\begin{equation}
K\circ \varphi_W=\deg \circ \alpha_W
\label{IKdH}
\end{equation}
in view of the definition of $\varphi_W$. This implies that one has
\begin{equation}
K\circ L_h=L_h\circ K
\label{KI}
\end{equation}
and
\begin{align*}
(K\circ d + d\circ K) \circ\alpha_W &=\alpha_W\circ \deg=\deg \circ \alpha_W
\\
(K\circ d + d\circ K)\circ d \circ\alpha_W &= d\circ \alpha_W\circ \deg=(\deg-I)\circ d \circ \alpha_W
\end{align*}
on $C(\calh)$ which implies $H^n(W(\calh))=0$ for $n\geq 1$ and by using (\ref{KI}), $H^n_I(W(\calh))=0$  for $n\geq 1$. On the other hand that $H^0(W(\calh))=H^0_I(W(\calh))=\mathbb K$ is obvious. $\qed$

The computation of the basic cohomology is more involved. This is connected with the fact that the structures of the horizontal subalgebra $W_H(\calh)$ and of the basic differential subalgebra $W_B(\calh)$ of $W(\calh)$ are more complicated than in the classical case  of the Weil algebra of a Lie algebra. We postpone the analysis of these points to the future Part II.

\section{Further comments}
This paper is the first part of a work on the noncommutative generalization of the notion of Cartan operation and of the Weil algebra. In this first part we have set up the general formulation of this noncommutative version. Part II will be devoted to the description in this context of the noncommutative version of the Weil homomorphism and of the noncommutative version of the Cartan map.  

Let us now explain why we did not mention the axiom $(i_X)^2=0$ (for $X\in \fracg)$ for the operation of a Lie algebra $\fracg$ in a graded differential algebra. Firstly, $(i_X)^2$ is (for $X\in\fracg$) a derivation of degree -2 of $\Omega$ which implies that it vanishes on the graded subalgebra of $\Omega$ generated by the elements of degrees 0 and 1 (remembering that $\Omega$ is positively graded by assuption).
Secondly the axiom \ref{axop} (Cartan relation) implies that one has
\[
[(i_X)^2,d]=0 \qquad \text{for} \quad X\in \fracg \, .
\]

Thus $(i_X)^2=0$ on the graded differential subalgebra $\Omega_{(1)}$ of $\Omega$ generated (as graded differential algebra) by the elements of degrees 0 and 1. In all cases of interest one has $\Omega_{(1)}=\Omega$, that is $\Omega$ is generated as graded differential algebra by $\Omega^0 \oplus \Omega^1$ which implies $(i_X)^2=0$ for $X\in \fracg$. Thus one does not need the axiom $(i_X)^2=0$ which plays no role otherwise.  

Finally let us say some words on the relation with the construction of 
\cite{cra:2002}. In the interesting paper \cite{cra:2002} there is a definition of the Weil algebra of a coalgebra which leads of course to a definition of a Weil algebra of a Hopf algebra $\calh$. However the corresponding Weil algebra is generated by $\calh$ instead of $\calh^\ast$ as our $W(\calh)$. Thus in spite of some similarities, it is a different object which is considered in \cite{cra:2002}. In particular our correspondence $\calh\mapsto W(\calh)$ has the same variance as the classical correspondence $\fracg\mapsto W(\fracg)$ $(W(\fracg)$ is generated by $\fracg^\ast$) while the Weil algebra of \cite{cra:2002} has an opposite variance.  In fact, the dual ``comodule algebras" approach is interesting and convenient for some purpose but is inappropriate for the generalization of the formulation of Henri Cartan in terms of operations and algebraic connections. Nevertheless, for the noncommutative Weil homomorphism, we shall use results of \cite{cra:2002} as well as those of \cite{ac-mos:1998}, \cite{ac-mos:1999} and \cite{ac-mos:2000}. \\ 

\noindent {\bf Acknowledgements.}
It is a pleasure to thank Jim Stasheff and Robert Coquereaux for their kind interest and advices. We are grateful to the referees for their constructive criticisms. 

\section*{Appendices}
\appendix

\section{The Hopf superalgebra formulation of operations}

Let $\calh$ be a Hopf algebra. Then a $\calh$-{\sl algebra} is an algebra $\cala$ endowed with an action $\calh\otimes \cala\rightarrow \cala$ of $\calh$, $h\otimes a\mapsto ha$ satisfying (cf. \cite{cra:2002}, \cite{ac-mos:1999})
\begin{equation}
h(ga)=(hg)a,\>\>\>\ \forall h, g\in \calh,\>\> a\in \cala
\label{pac}
\end{equation}
and 
\begin{equation}
h(ab)=\sum\nolimits_j(h^{(1)}_j a)(h^{(2)}_j b),\>\>\>\> \forall h\in \calh,\> a,b\in \cala \, .
\label{cac}
\end{equation} 

There is an obvious graded (super) version of the above notion: Let $\calh$ be a $\mathbb Z$-graded (super) Hopf algebra, then a {\sl graded $\calh$-algebra} is a $\mathbb Z$-graded algebra $\cala$ with a {\sl homogeneous} action $\calh\otimes \cala\rightarrow \cala$ satisfying (\ref{pac}) and (\ref{cac}). \\

Let $\calh$ be an ordinary (non graded) Hopf algebra. To $\calh$ one associates a $\mathbb Z$-graded Hopf algebra $\widehat\calh$ in the following manner. In degree 1, $\widehat\calh_1$ is 1-dimensional generated by an element $\delta$ with properties 
\[
\delta^2=0,\>\>\> \Delta\delta=\delta\otimes \bbbone +\bbbone \otimes \delta,\>\>\> \varepsilon(\delta)=0,\>\>\> S(\delta)=-\delta
\]
where $\Delta$ is the coproduct, $\varepsilon$ is the counit  and $S$ is the antipode of $\widehat\calh$. In degree 0, $\widehat\calh_0$ is isomorphic to $\calh$ as Hopf algebra and we denote by
\[
\Lambda:\calh\mapsto \widehat\calh_0,\>\>\> h\mapsto \Lambda_h
\]
 the corresponding isomorphism. The relations with $\delta$ are
 \[
 [\delta,\Lambda_h]=0  
 \]
for $h\in \calh $.
  In degree -1, $\widehat\calh_{-1}$ is isomorphic as vector space to the quotient $\calh/\mathbb K\bbbone$ and we denote by 
 \[
 y:\calh\mapsto \widehat\calh_{-1},\>\>\> h\mapsto y_h
 \]
  the corresponding linear mapping which vanishes on $\bbbone\in \calh$ (i.e. $y_{\bbbone}=0$). The relations, the coproduct and the antipode of $y_h$ are given by 
\begin{align*}
  y_g\Lambda_h &= \sum\nolimits_j\Lambda_{h^{(1)}_j}y_{\ad(h^{(2)}_j)g}   \\
 {[}\delta,y_h] &=  \Lambda_h-\varepsilon(h)\bbbone  \\
 \Delta y_h &=  \sum\nolimits_jy_{h^{(1)}_j}\otimes \Lambda_{h^{(2)}_j}+\bbbone\otimes y_h  \\
 S(y_h) &= -\sum\nolimits_j y_{h^{(1)}_j} \Lambda_{S(h^{(2)}_j)}
\end{align*}
for $g,h\in \calh$,  while one has $\varepsilon(y_h)=0$. The graded Hopf algebra $\widehat\calh$ is generated by $\widehat\calh_1\oplus \widehat\calh_0\oplus \widehat\calh_{-1}$ and there are no other relations in $\widehat\calh$.  

Let us next consider an operation of $\calh$ in a graded differential algebra $(\Omega,d)$. Then $\delta\mapsto d, \Lambda_h\mapsto L_h$ and $y_h\mapsto i_h$ define a structure of graded $\widehat\calh$-algebra on $\Omega$. It is clear that this correspondence allows to identify the notion of $\calh$-operation with the notion of positively graded $\widehat\calh$-algebra.  
This is the counterpart in our noncommutative context of the classical graded Lie superalgebra formulation of the operations of Lie algebras \cite{mat-qui:1986}, \cite{gui-ste:1999}, \cite{ale-mei:1999}.

\section{More on algebraic connections in $\fracg$-operations}

Let $\Omega$ be a graded commutative differential algebra which is endowed with an operation of the Lie algebra $\fracg$. 
Then, to give a connection $\alpha$ in $\Omega$ as in \S\ref{acgp} is the same as to give an element $A$ of $\fracg\otimes \Omega^1$ such that
\begin{equation}
i_X(A)=X \qquad 
\text{and} \qquad
L_X(A)=\ad(X)A
\label{AclL}
\end{equation}
where here $i_X$ is $I_\fracg\otimes i_X$, $L_X$ is $I_\fracg\otimes L_X$ and $\ad(X)$ is $\ad(X)\otimes I_{\Omega^1}$. The mapping $\alpha$ being then $\theta\mapsto (\theta\otimes I_{\Omega^1})A$.  

In the same vein, the curvature $\varphi$  is the same as the element $F$ of $\fracg\otimes \Omega^2$ defined by
\begin{equation}
F=dA+\tfrac{1}{2}[A,A]
\label{AFcurv}
\end{equation}
where now $d=I_\fracg\otimes d$ and $[\>,\>]$ is the graded commutator in $\fracg\otimes \Omega$, the mapping $\varphi$ being then $\theta\mapsto (\theta\otimes I_{\Omega^2})F$.  One has in view of both 
(\ref{InvC})
\begin{equation}
i_X(F)=0
\qquad \text{and} \qquad
L_X(F)=\ad(X)F
\label{InvCF}
\end{equation}
for $X\in \fracg$ with obvious notations.  Bianchi identity is easily established: $d F + [A, F]=0.$

\section{The case $\dim(\calh)<\infty$ and the case $\dim(\fracg)=\infty$} \label{cases}

If $\dim(\calh)<\infty$, then $C(\calh)$ is the tensor algebra $T(\calh^\ast)$ and a connection $\alpha$ is completely specified by its restriction to $\calh^\ast$. Thus in this case, a connection is simply a linear mapping
\begin{equation}
\alpha:\calh^\ast\rightarrow \Omega^1
\label{fCo}
\end{equation}
satisfying {both relations in 
(\ref{dCo2}).  

Since in this case one has the isomorphism 
\[
\Hom(\calh^\ast,\Omega^1)\simeq \calh\otimes \Omega^1
\]
of vector spaces, one can also say that a connection on the operation of $\calh$ in $\Omega$ is an element $A$ of degree 1 of $\calh\otimes\Omega$, that is 
\begin{equation}
A\in \calh\otimes \Omega^1
\label{Co0}
\end{equation}
such that, for any $h\in \calh$, 
\begin{equation}
i_h(A)=h-\varepsilon(h)\bbbone
\qquad \text{and} \qquad
L_h(A)=\text{ad}(h)A \, .
\label{Co2}
\end{equation}

Given such a connection $A$, {\sl its curvature} $F=F(A)$ is the element
\begin{equation}
F=dA+A^2-(\varepsilon A + A\varepsilon)=d(A-\varepsilon)+(A-\varepsilon)^2
\label{curv}
\end{equation}
of $\calh\otimes \Omega^2$ corresponding to $\varphi\in\Hom(\calh^\ast,\Omega^2)$.  
The curvature $F$ of $A$ satisfies
\begin{equation}
i_h(F)=0
\qquad \text{and} \qquad
L_h(F)=\mathrm{ad}(h)F
\label{equiF}
\end{equation}
for any $h\in \calh$. This can be checked directly by using (\ref{curv}), 
both (\ref{Co2}) and the definition of $\calh$-operations. Applying $d$ to (\ref{curv}) implies 
\begin{equation}
dF+(A-\varepsilon) F-F(A-\varepsilon)=0
\label{BianchF}
\end{equation}
which is the Bianchi identity in the present context.  

One sees that the case where $\dim(\calh)<\infty$ looks formally close to the classical theory summarized in Section~\ref{cathe}. However, it is worth noticing here that in Section~\ref{cathe} as well as in the original references \cite{car:1951a} and \cite{car:1951b}, it is implicitely assumed that the Lie algebra $\fracg$ is finite-dimensional. If one wants to extend the classical theory (where all graded algebras are graded commutative) to the case where $\dim(\fracg)=\infty$, then one must change accordingly the definition of algebraic connections. Indeed, the transposed of the Lie bracket is now a linear mapping of $\fracg^\ast$ into the space $C^2_\wedge(\fracg)$ of antisymmetric bilinear forms on $\fracg$ and therefore one has to replace $\wedge\fracg^\ast$ by the graded commutative algebra
\[
C_\wedge(\fracg)=\oplus_n C^n_\wedge (\fracg)
\]
where $C^n_\wedge(\fracg)$ is the vector space of completely antisymmetric $n$-linear forms on $\fracg$ and one endows $C_\wedge(\fracg)$ with the Chevalley-Eilenberg differential $d$ of cochains on $\fracg$ with values in the trivial representation (in $\mathbb K$). Then, given an operation $i$ of $\fracg$ in the graded commutative differential algebra $\Omega=\oplus_n \Omega^n$, an algebraic connection on $\Omega$ should be defined as a homomorphism of graded commutative algebras
\[
\alpha:C_\wedge(\fracg)\rightarrow \Omega
\]
satisfying conditions similar to 
(\ref{GCo2}) with curvature $\varphi$ again defined by
\[
\varphi=d\circ \alpha-\alpha\circ d
\]
with obvious notations, etc.  
Thus for $\dim(\fracg)=\infty$, all definitions and formulas look completely similar to the one of \S \ref{alco}, except that we are in a graded commutative context with the classical notion of operation of the Lie algebra $\fracg$.  

The definition of the Weil algebra $W(\fracg)$ of an infinite-dimensional Lie algebra $\fracg$ must be modified accordingly and looks then closer to the definition of the Weil algebra $W(\calh)$ of a Hopf algebra $\calh$ defined in Section~\ref{uhowc}.


\begin{thebibliography}{10}

\bibitem{ale-mei:1999}
A.~Alekseev and E.~Meinrenken.
\newblock The non-commutative {W}eil algebra.
\newblock {\em Inv. Math.}, 139:135--172, 2000.

\bibitem{asc-cas:1993}
P.~Aschieri and L.~Castellani.
\newblock An introduction to noncommutative differential geometry on quantum
  groups.
\newblock {\em Int. J. Mod. Phys.}, A8:1667--1706, 1993.

\bibitem{bou:1970}
N.~Bourbaki.
\newblock {\em Alg{\`e}bre {I}, {C}hapitre {III}}.
\newblock Hermann, 1970.

\bibitem{bra-lan:2011}
S.~Brain and G.~Landi.
\newblock Differential and twistor geometry of the quantum {H}opf fibration.
\newblock ArXiv:1103.0419.

\bibitem{brz-maj:1993}
T.~Brzezi\'nski and S.~Majid.
\newblock Quantum group gauge theory on quantum spaces.
\newblock {\em Commun. Math. Phys.}, 157:591--638, 1993.

\bibitem{brz-maj:1995}
T.~Brzezi\'nski and S.~Majid.
\newblock Erratum.
\newblock {\em Commun. Math. Phys.}, 167:235, 1995.

\bibitem{car:1951a}
H.~Cartan.
\newblock \phantom{a} {N}otion d'alg{\`e}bre diff{\'e}rentielle; application
  aux groupes de {L}ie et aux vari{\'e}t{\'e}s o{\`u} op{\`e}re un groupe de
  {L}ie.
\newblock In {\em Colloque de Topologie, Bruxelles, 1950}, pages 15--27.
  Masson, 1951.

\bibitem{car:1951b}
H.~Cartan.
\newblock \phantom{b} {L}a trangression dans un groupe de {L}ie et dans un
  espace fibr{\'e} principal.
\newblock In {\em Colloque de Topologie, Bruxelles, 1950}, pages 57--71.
  Masson, 1951.

\bibitem{car-eil:1973}
H.~Cartan and S.~Eilenberg.
\newblock {\em Homological algebra}.
\newblock Princeton University Press, 1973.

\bibitem{ac:1982}
A.~Connes.
\newblock Noncommutative differential geometry. {P}art {I}: The {C}hern
  character in ${K}$-homology; {P}art {II}: de {R}ham homology and
  noncommutative algebra.
\newblock {\em Preprints IHES M/82/53; M83/19}, 1982/83.

\bibitem{ac:1994}
A.~Connes.
\newblock {\em Non-commutative geometry}.
\newblock Academic Press, 1994.

\bibitem{ac-mdv:2002a}
A.~Connes and M.~Dubois-Violette.
\newblock Noncommutative finite-dimensional manifolds. {I}. {S}pherical
  manifolds and related examples.
\newblock {\em Commun. Math. Phys.}, 230:539--579, 2002.

\bibitem{ac-lan:2001}
A.~Connes and G.~Landi.
\newblock Noncommutative manifolds, the instanton algebra and isospectral
  deformations.
\newblock {\em Commun. Math. Phys.}, 221:141--159, 2001.

\bibitem{ac-mos:1998}
A.~Connes and H.~Moscovici.
\newblock Hopf algebras, cyclic cohomology and the transverse index theorem.
\newblock {\em Commun. Math. Phys.}, 198:199--246, 1998.

\bibitem{ac-mos:1999}
A.~Connes and H.~Moscovici.
\newblock Cyclic cohomology and {H}opf algebras.
\newblock {\em Lett. Math. Phys.}, 48:97--108, 1999.

\bibitem{ac-mos:2000}
A.~Connes and H.~Moscovici.
\newblock Cyclic cohomology and {H}opf algebra symmetry.
\newblock {\em Lett. Math. Phys.}, 53:1--28, 2000.

\bibitem{coq-kas:1989}
R.~Coquereaux and D.~Kastler.
\newblock Remarks on the differential envelopes of associative algebras.
\newblock {\em Pacific J. Math.}, 137:245--263, 1989.

\bibitem{cra:2002}
M.~Crainic.
\newblock Cyclic cohomology of {H}opf algebras.
\newblock {\em Journal of Pure and Applied Algebra}, 166:29--66, 2002.

\bibitem{mdv:1987b}
M.~Dubois-Violette.
\newblock The {W}eil-{B.R.S.} algebra of a {L}ie algebra and the anomalous
  terms in gauge theory.
\newblock {\em J. Geom. Phys.}, 3:525--565, 1987.

\bibitem{mdv:1988}
M.~Dubois-Violette.
\newblock D\'erivations et calcul diff\'erentiel non-commutatif.
\newblock {\em C.R. Acad. Sci. Paris}, 307:403--408, 1988.

\bibitem{mdv:2001}
M.~Dubois-Violette.
\newblock Lectures on graded differential algebras and noncommutative geometry.
\newblock In Y.~Maeda and al., editors, {\em Noncommutative Differential
  Geometry and Its Applications to Physics}, pages 245--306. Shonan, Japan,
  1999, Kluwer Academic Publishers, 2001.

\bibitem{mdv-hen-tal:1992}
M.~Dubois-Violette, M.~Henneaux, M.~Talon, and C.M. Viallet.
\newblock General solution of the consistency equation.
\newblock {\em Phys. Lett.}, B289:361--367, 1992.

\bibitem{mdv-lan:2011}
M.~Dubois-Violette and G.~Landi.
\newblock Lie prealgebras.
\newblock In A.~Connes and al., editors, {\em Noncommutative {G}eometry and
  {G}lobal {A}nalysis}, volume 546 of {\em Contemporary Mathematics}, pages
  115--135. American Mathematical Society, 2011.

\bibitem{mdv-tal-via:1985b}
M.~Dubois-Violette, M.~Talon, and C.M. Viallet.
\newblock {B.R.S.} algebras. analysis of consistency equations in gauge theory.
\newblock {\em Commun. Math. Phys}, 102:105 --122, 1985.

\bibitem{gre-hal-van:1976}
W.~Greub, S.~Halperin, and R.~Vanstone.
\newblock {\em Connections, curvature, and cohomology}, volume III.
\newblock Academic Press, 1976.

\bibitem{gui-ste:1999}
V.~Guillemin and S.~Sternberg.
\newblock {\em Supersymmetry and equivariant de {R}ham theory}.
\newblock Springer Verlag, 1999.

\bibitem{kam-ton:1975}
F.W. Kamber and P~Tondeur.
\newblock {\em Foliated bundles and characteristic classes}, volume 993 of {\em
  Lecture Notes in Mathematics}.
\newblock Springer Verlag, 1975.

\bibitem{kar:1983}
M.~Karoubi.
\newblock Homologie cyclique des groupes et alg{\`e}bres.
\newblock {\em C.R. Acad. Sci. Paris, S{\'e}rie {I}}, 297:381--384, 1983.

\bibitem{kar:1987}
M.~Karoubi.
\newblock {\em Homologie cyclique et {K}-th{\'e}orie}, volume 149 of {\em
  Ast\'erique}.
\newblock Soci{\'e}t\'e Math\'ematique de France, 1987.

\bibitem{kli-sch:1997}
A.~Klimyk and K.~Schm{\"u}dgen.
\newblock {\em Quantum groups and their representations}.
\newblock Springer, 1997.

\bibitem{kob-nom:1963}
M.~Kobayashi and K.~Nomizu.
\newblock {\em Foundations of differential geometry}, volume~1.
\newblock Interscience publishers, 1963.

\bibitem{lan-sui:2005}
G.~Landi and W.D. van Suijlekom.
\newblock Principal fibrations from noncommutative spheres.
\newblock {\em Commun. Math. Phys.}, 260:203--225, 2005.

\bibitem{lan-sui:2008}
G.~Landi and W.D. van Suijlekom.
\newblock Noncommutative instantons in {T}ehran.
\newblock In M.~Khalkhali and al., editors, {\em An invitation in
  noncommutative geometry}, pages 275--353, Hackensack, N.J., 2008. World
  Scientific Pub.

\bibitem{lod:1992}
J.L. Loday.
\newblock {\em Cyclic homology}, volume 301 of {\em Grundlehren der
  mathematischen Wissenschaften}.
\newblock Springer Verlag, New York, 1992.

\bibitem{mat-qui:1986}
V.~Mathai and D.~Quillen.
\newblock Thom classes superconnections and equivariant differential forms.
\newblock {\em Topology}, 25:85--106, 1986.

\bibitem{sch:2004}
P.~Schauenburg.
\newblock Hopf-{G}alois and bi-{G}alois extensions.
\newblock In G.~Janelidze and al., editors, {\em Galois theory, {H}opf algebras
  and semiabelian categories}, volume~43 of {\em Fields Institute
  Communications}, pages 469--515. American Mathematical Society, 2004.

\bibitem{sch-wat:1994}
P.~Schupp and P.~Watts.
\newblock Universal and generalized {C}artan calculus on {H}opf algebras.
\newblock ArXiv:hep-th/9402134.

\bibitem{sul:1977}
D.~Sullivan.
\newblock Infinitesimal computations in topology.
\newblock {\em Pub. IHES}, 47:269--331, 1977.

\bibitem{wor:1987b}
S.L. Woronowicz.
\newblock Twisted {SU}(2) group. an example of noncommutative differential
  calculus.
\newblock {\em Publ. RIMS, Kyoto Univ.}, 23:117--181, 1987.

\bibitem{wor:1989}
S.L. Woronowicz.
\newblock Differential calculus on compact matrix pseudogroups (quantum
  groups).
\newblock {\em Commun. Math. Phys.}, 122:129--170, 1989.

\end{thebibliography}

\end{document}